\documentclass[a4,12pt]{amsart}
\usepackage{amsfonts,amssymb,amscd,amsmath,enumerate,verbatim,calc}
\newtheorem{Theorem}{Theorem}[section]
\newtheorem{Lemma}[Theorem]{Lemma}
\newtheorem{Corollary}[Theorem]{Corollary}

\theoremstyle{definition}
\newtheorem{Remark}[Theorem]{Remark}

\newtheorem{Example}[Theorem]{Example}

\newtheorem{Definition}[Theorem]{Definition}

\def\RMN#1{\uppercase\expandafter{\romannumeral#1}}

\def\uop#1{\mathop{#1}\nolimits}
\def\standop#1{\uop{\mathrm{#1}}}
\def\difstop#1#2{\expandafter\def\csname #1\endcsname{\standop{#2}}}
\def\defstop#1{\difstop{#1}{#1}}
\difstop{inn}{Int}
\difstop{aut}{Aut}
\difstop{kaku}{Ker}
\difstop{zou}{Im}

\DeclareMathOperator{\spec}{Spec}
\let\Spec\spec
\DeclareMathOperator{\proj}{Proj}

\makeatletter
  
  \@addtoreset{equation}{section}
\makeatother

\renewcommand{\qed}{{\unskip\nobreak\hfil\penalty50\quad\null\nobreak\hfil{\bf
q.e.d.}\parfillskip0pt\finalhyphendemerits0\par\medskip}}
\defstop{Cl}
\defstop{Lqc}
\defstop{Sing}
\defstop{Qch}
\defstop{Coh}
\defstop{Cox}
\defstop{Hom}
\def\N{\mathcal{N}}
\def\M{\mathcal{M}}
\def\A{\mathcal{A}}
\def\O{\mathcal{O}}
\def\E{\mathcal{E}}
\def\fm{\mathfrak{m}}
\let\Bbb\mathbb
\let\frak\mathfrak
\def\uH{\mathop{\mbox{\underline{$H$}}}\nolimits}
\def\uHom{\mathop{\mbox{\underline{$\mathrm{Hom}$}}}\nolimits}
\def\uExt{\mathop{\mbox{\underline{$\mathrm{Ext}$}}}\nolimits}

\begin{document}

\title{The canonical module of a Cox ring}
\author{Mitsuyasu Hashimoto and Kazuhiko Kurano}
\maketitle

\begin{center}
{\small Dedicated to Professor Jun-ichi Nishimura on the occasion of his sixtieth birthday}
\end{center}

\begin{abstract}
In this paper, we shall describe the graded canonical module of a Noetherian multi-section ring
of a normal projective variety.
In particular, in the case of the Cox ring, we prove that the graded canonical module
is a graded free module of rank one with the shift of degree $K_X$.
We shall give two kinds of proofs.
The first one utilizes the equivariant twisted inverse functor developed by the first author.
The second proof is down-to-earth, that avoids the twisted inverse functor,
but some additional assumptions are required in this proof.
\end{abstract}

\section{Introduction}\label{Intro}

D.~Cox studied the total coordinate ring (Cox ring) of toric varieties in \cite{cox-hom}.
Cox rings of normal projective varieties have become so important and 
interesting objects that many mathematicians
try to prove (in)finite generation or to study their ring-theoretic properties (generators, relations, 
syzygies, homological properties, etc.). 

In this paper, we shall describe the graded canonical module of a Noetherian multi-section ring
of a normal projective variety in Theorem~\ref{main}.
Using this result, we shall give a necessary and sufficient condition for the canonical module 
to be a free module (see Corollary~\ref{k}).
Since the Cox ring is a unique factorization domain as in \cite{EKW},
its graded canonical module is a free module.
We prove that the graded canonical module of the Cox ring of a normal projective variety $X$
is a graded free module of rank one with the shift of degree $K_X$ (see Corollary~\ref{cor}).
Many mathematicians study syzygies of the Cox ring.
One of the purposes of our study is to contribute to their study.

We shall give two kinds of proofs to Theorem~\ref{main} in this paper.
In Section~\ref{proof}, we prove Theorem~\ref{main} using 
the equivariant twisted inverse functor developed in \cite{LH}.
The second proof is given in Section~\ref{another}.
The second proof is down-to-earth, that avoids the twisted inverse functor,
but some additional assumptions are required in this proof.

From now on, we shall give precise definition, and state our main theorem and corollaries.

\begin{Definition}
\begin{rm}
Let $X$ be a $d$-dimensional normal projective variety over a field $k$.
Let $D_1$, \ldots, $D_s$ be Weil divisors on $X$.
We define a ring $R(X; D_1, \ldots, D_s)$ to be
\[
\bigoplus_{(n_1, \ldots, n_s) \in {\mathbb Z}^s} H^0(X, {\mathcal O}_X(\sum_i n_iD_i))t_1^{n_1}\cdots t_s^{n_s}
\subset k(X)[t_1^{\pm 1}, \ldots, t_s^{\pm 1}] .
\]
For a Weil divisor $F$ on $X$, we set
\[
M_F = 
\bigoplus_{(n_1, \ldots, n_s) \in {\mathbb Z}^s} H^0(X, {\mathcal O}_X(\sum_i n_iD_i +F))
t_1^{n_1}\cdots t_s^{n_s}
\subset k(X)[t_1^{\pm 1}, \ldots, t_s^{\pm 1}] ,
\]
that is, $M_F$ is a ${\mathbb Z}^s$-graded $R(X; D_1, \ldots, D_s)$-module such that
\[
[M_F]_{(n_1, \ldots, n_s)} = H^0(X, {\mathcal O}_X(\sum_i n_iD_i +F))
t_1^{n_1}\cdots t_s^{n_s} .
\]
\end{rm}
\end{Definition}

In this paper, for a normal variety $X$, we denote by 
$\Cl(X)$ the class group of $X$, and
for a Weil divisor $F$ on $X$,
we denote by $\overline{F}$ the class of $F$ 
in ${\rm Cl}(X)$.

In the case where ${\rm Cl}(X)$ is freely generated by $\overline{D_1}$, \ldots, $\overline{D_s}$,
the ring $R(X; D_1, \ldots, D_s)$ is called the {\em Cox ring} of $X$
and denoted by ${\rm Cox}(X)$.
Recall that ${\rm Cox}(X)$ is uniquely determined by $X$ up to isomorphisms,
that is, independent of the choice of  $D_1, \ldots, D_s$.

\begin{Theorem}\label{main}
Let $X$ be a normal projective variety over a field $k$ such that $H^0(X, {\mathcal O}_X) = k$.
Assume that $D_1$, \ldots, $D_s$ are Weil divisors on $X$ which satisfy the following three conditions:
\begin{enumerate}
\item
$\overline{D_1}$, \ldots, $\overline{D_s}$ are linearly independent over ${\mathbb Z}$ in the divisor class group ${\rm Cl}(X)$.
\item
The lattice ${\mathbb Z}D_1 + \cdots + {\mathbb Z}D_s$ contains an ample Cartier divisor.
\item
The ring $R(X; D_1, \ldots, D_s)$ is Noetherian.
\end{enumerate}

Then, $R(X; D_1, \ldots, D_s)$ is a local ${\mathbb Z}^s$-graded $k$-domain 
%%% modified 2011/03/11
{\rm(}see Definition~\ref{local.def}{\rm)},
and we have an isomorphism
\[
\omega_{R(X; D_1, \ldots, D_s)}
\simeq M_{K_X}
\]
of ${\mathbb Z}^s$-graded $R(X; D_1, \ldots, D_s)$-modules.
\end{Theorem}

Let $X$ be a normal projective variety  and $D_1$, \ldots, $D_s$ be Weil divisors  on $X$.
Then, $R(X; D_1, \ldots, D_s)$ is a Krull domain as in \cite{EKW}.
(In particular, if $R(X; D_1, \ldots, D_s)$ is Noetherian, then it is a Noetherian
integrally closed domain.)
Therefore, the divisor class group of $R(X; D_1, \ldots, D_s)$ can be defined.
%If the field of quotients of $R(X; D_1, \ldots, D_s)$ is equal to 
%$k(X)(t_1, \ldots, t_s)$, 
%then we have a surjective map 
%\[
%p : {\rm Cl}(X) \longrightarrow {\rm Cl}(R(X; D_1, \ldots, D_s))
%\]
%defined by 
%\begin{equation}\label{p}
%p(\overline{F}) = \overline{{M_F}^{**}} ,
%\end{equation}
%where the right hand side stands for the isomorphism class of ${M_F}^{**}$
%and $[ \ \ ]^*$ is the $R(X; D_1, \ldots, D_s)$-dual (see \cite{EKW}).
%Here, remark that the map $\varphi$ in Theorem~1.1 in \cite{EKW} equals to $-p$.
%If we do not assume the condition (2) in Theorem~\ref{main},
%then the height of a prime ideal $M_{-V}$ in (\ref{P_V}) can be bigger than one (p628 in \cite{EKW}).
%Therefore, we have to take the double dual in (\ref{p}).
%
%Further, if  the assumption (2) in Theorem~\ref{main} is satisfied, then
%\begin{equation}\label{ekw}
%0 \longrightarrow {\mathbb Z}\overline{D_1} + \cdots + {\mathbb Z}\overline{D_s}
%\longrightarrow {\rm Cl}(X)
%\stackrel{p}{\longrightarrow} {\rm Cl}(R(X; D_1, \ldots, D_s)) \longrightarrow 0
%\end{equation}
%is exact as in \cite{EKW}.
%(If the assumption (2) in Theorem~\ref{main} is satisfied, 
%then the field of quotients of $R(X; D_1, \ldots, D_s)$ coincides with 
%$k(X)(t_1, \ldots, t_s)$,
%and $M_F$ is a reflexive $R(X; D_1, \ldots, D_s)$-module 
%for each Weil divisor $F$ on $X$ as in Lemma~\ref{lem2}.)
If (2) in Theorem~\ref{main} is satisfied,
then we have an exact sequence
\begin{equation}\label{ekw}
0 \longrightarrow {\mathbb Z}\overline{D_1} + \cdots + {\mathbb Z}\overline{D_s}
\longrightarrow {\rm Cl}(X)
\stackrel{p}{\longrightarrow} {\rm Cl}(R(X; D_1, \ldots, D_s)) \longrightarrow 0
\end{equation}
defined by $p(\overline{F}) = \overline{M_F}$ for each Weil divisor $F$,
where $\overline{M_F}$ stands for the isomorphism class
containing $M_F$. 
We shall construct this map $p$ explicitly in Remark~\ref{divisorial}.

Theorem~\ref{main} says that 
\[
p(\overline{K_X}) = \overline{\omega_{R(X; D_1, \ldots, D_s)}}
\]
if (1), (2), and (3) in Theorem~\ref{main} are satisfied.
In particular, we have the following corollary:

\begin{Corollary}\label{k}
Suppose that the assumptions in Theorem~\ref{main} are satisfied.

Then, $\omega_{R(X; D_1, \ldots, D_s)}$ is a free $R(X; D_1, \ldots, D_s)$-module if and only if 
\[
\overline{K_X} \in {\mathbb Z}\overline{D_1} + \cdots + {\mathbb Z}\overline{D_s} 
\]
in ${\rm Cl}(X)$.
\end{Corollary}

We shall give some examples (Example~\ref{ex1}, Example~\ref{ex2}) in Section~\ref{sectionex}.

\begin{Remark}
\begin{rm}
Suppose that the assumptions in Theorem~\ref{main} are satisfied.
Put $d = \dim X$.
Let $F_1$, \ldots, $F_r$ be linearly independent elements 
in ${\mathbb Z}D_1 + \cdots + {\mathbb Z}D_s$.
Put
\[
L_1 = {\mathbb Z}D_1 + \cdots + {\mathbb Z}D_s
\supset L_2 = {\mathbb Z}F_1 + \cdots + {\mathbb Z}F_r .
\]
We think that $R(X; D_1, \ldots, D_s)$ is graded by $L_1$.
We have
\[
R(X; F_1, \ldots, F_r) = R(X; D_1, \ldots, D_s)|_{L_2} .
\]
By Theorem~\ref{main},
\begin{equation}\label{seigen}
\omega_{R(X; F_1, \ldots, F_r)} = \omega_{R(X; D_1, \ldots, D_s)}|_{L_2} 
\end{equation}
if the lattice $L_2$ contains an ample Cartier divisor.

If we take the graded dual on the both sides of (\ref{seigen}),
we know
\begin{equation}\label{seigencohomology}
H^{d+r}_{{\mathfrak m}_2}(R(X; F_1, \ldots, F_r)) = 
H^{d+s}_{{\mathfrak m}_1}(R(X; D_1, \ldots, D_s))|_{L_2},
\end{equation}
where ${\mathfrak m}_1$ (resp.\ ${\mathfrak m}_2$) is the unique maximal
homogeneous ideal of $R(X; D_1, \ldots, D_s)$ (resp.\ $R(X; F_1, \ldots, F_r)$).
Here, remark that the dimension of $R(X; D_1, \ldots, D_s)$ (resp.\ $R(X; F_1, \ldots, F_r)$)
is equal to $d+s$ (resp.\ $d+r$). 
We shall obtain the equalities
\begin{eqnarray*}
H^{d+s}_{{\mathfrak m}_1}(R(X; D_1, \ldots, D_s)) & = & H^{d+1}_{S_+}(R(X; D_1, \ldots, D_s))  \\
H^{d+r}_{{\mathfrak m}_2}(R(X; F_1, \ldots, F_r)) & = & H^{d+1}_{S_+}(R(X; F_1, \ldots, F_r)) 
\end{eqnarray*}
in Remark~\ref{4.3}.
The equality (\ref{seigencohomology}) 
also follows from the above equalities.

In Example~\ref{ex3}, we shall give an example that the equality~(\ref{seigen}) is not satisfied 
if we remove the assumption that $L_2$ contains an ample Cartier divisor.
\end{rm}
\end{Remark}

Suppose that $X$ is a normal projective variety with a 
finitely generated free divisor class group.
We may think that the Cox ring is graded by ${\rm Cl}(X)$.
By the exact sequence (\ref{ekw}),
${\rm Cox}(X)$ is a unique factorization domain as in \cite{EKW}.
Therefore the canonical module $\omega_{{\rm Cox}(X)}$ is a ${\rm Cl}(X)$-graded free 
${\rm Cox}(X)$-module of rank one if ${\rm Cox}(X)$ is a Noetherian ring.
By Theorem~\ref{main}, we can determine  the degree of the homogeneous generator 
of  $\omega_{{\rm Cox}(X)}$ as follows.

\begin{Corollary}\label{cor}
Let $X$ be a normal projective variety over a field.
Assume that ${\rm Cl}(X)$ is a finitely generated free abelian group
and the Cox ring of $X$ is Noetherian.

Then, 
%we have an isomorphism
%\[
%\omega_{{\rm Cox}(X)}
% \simeq \bigoplus_{(n_1, \ldots, n_s) \in {\mathbb Z}^s} H^0(X, {\mathcal O}_X(\sum_i n_iD_i + K_X)) 
%\]
%of  ${\rm Cl}(X)$-graded free ${\rm Cox}(X)$-modules,
%that is, 
the canonical module of the Cox ring is a rank one free ${\rm Cl}(X)$-graded module whose 
generator is of degree $-\overline{K_X} \in {\rm Cl}(X)$.
\end{Corollary}

Suppose that $D_1$, \ldots, $D_s$ are ${\mathbb Q}$-divisors.
Assume (1), (2), and (3) in Theorem~\ref{main}.
Let $H_1$ be the set of all the closed subvarieties of $X$ of codimension one.
Set 
\[
D_i = \sum_{V\in H_1} \frac{p_{i, V}}{q_{i, V}}V,
\]
where $p_{i, V}$'s and $q_{i, V}$'s are integers such that $(p_{i, V}, q_{i, V}) =1$ and $q_{i, V} > 0$ for each $V$.
(If $p_{i, V} = 0$, then $q_{i, V} =1$.)
Then, as in Theorem~(2.8) in \cite{Wa}, we obtain the following corollary:

\begin{Corollary}
With the notation as above,
\[
\omega_{R(X; D_1, \ldots, D_s)}
\simeq \bigoplus_{(n_1, \ldots, n_s) \in {\mathbb Z}^s} H^0(X, {\mathcal O}_X(\sum_i n_iD_i + K_X 
+ \sum_{V\in H_1} \frac{q_V - 1}{q_V}V)) ,
\]
where $q_V$ is the least common multiple of $q_{1, V}$, \ldots, $q_{s, V}$ for each $V \in H_1$.
\end{Corollary}

\section{The canonical module of a local ${\mathbb Z}^s$-graded $k$-domain}\label{preliminary}

\begin{Definition}\label{local.def}
\begin{rm}
Let $k$ be a field.
$R$ is called a {\em local ${\mathbb Z}^s$-graded $k$-domain}
if the following conditions are satisfied:
\begin{itemize}
\item
$R = \bigoplus_{\underline{a} \in {\mathbb Z}} 
R_{\underline{a}}$ is a Noetherian ${\mathbb Z}^s$-graded domain.
\item
$R_{\underline{0}} = k$.
\item
Suppose that ${\frak m}$ is the ideal of $R$ generated by all the homogeneous elements of $R$
of degree different from $\underline{0}$. 
Then ${\frak m} \neq R$.
\end{itemize}
\end{rm}
\end{Definition}

Assume that $R$ is a local ${\mathbb Z}^s$-graded $k$-domain.
We remark that the ideal ${\frak m}$ as above is the unique maximal homogeneous ideal 
of  $R$.
Further, $R$ is of finite type over $k$.
Therefore, the height of ${\frak m}$ coincides with the dimension of $R$.

For a ${\mathbb Z}^s$-graded module $M$ over 
a local ${\mathbb Z}^s$-graded $k$-domain $R$,
${}^*{\rm Hom}_k(M,k)$ denotes the graded dual of $M$, that is,
${}^*{\rm Hom}_k(M,k)$ is a ${\mathbb Z}^s$-graded $R$-module such that
${}^*{\rm Hom}_k(M,k)_{\underline{a}} = {\rm Hom}_k(M_{-\underline{a}},k)$.

\begin{Definition}\label{canonical}
\begin{rm}
Let $R$ be a local ${\mathbb Z}^s$-graded $k$-domain with the maximal homogeneous ideal ${\frak m}$.
Then, 
\[
{}^*{\rm Hom}_k(H^{\dim R}_{{\frak m}}(R),k)
\]
is called the {\em canonical module} of $R$, and denoted by $\omega_R$.
Here, $H^{\dim R}_{{\frak m}}(R)$ is the $(\dim R)$-th local cohomology group and
it has a natural structure of a ${\mathbb Z}^s$-graded $R$-module.
\end{rm}
\end{Definition}

We emphasize
that $\omega_R$ has a structure of a ${\mathbb Z}^s$-graded $R$-module.
We refer the reader to \cite{GW2} and \cite{Ka} for the general theory of 
${\mathbb Z}^s$-graded rings.

\begin{Lemma}\label{local}
Let $X$ be a normal projective variety over a field $k$ such that $H^0(X, {\mathcal O}_X) = k$.
Assume that $D_1$, \ldots, $D_s$ are Weil divisors on $X$ which satisfy the assumptions 
(1) and (3) in Theorem~\ref{main}.

Then, $R(X; D_1, \ldots, D_s)$ is a local ${\mathbb Z}^s$-graded $k$-domain.
\end{Lemma}

\proof 
We denote the ring $R(X; D_1, \ldots, D_s)$ simply by $R$.

In order to prove that $R$ is a local ${\mathbb Z}^s$-graded $k$-domain,
we need to show that the ideal ${\frak m}$ is not equal to $R$, 
where ${\frak m}$ is the ideal of $R$ generated by all the homogeneous elements of $R$
of degree different from $\underline{0}$. 
Assume the contrary.
Then, there exists $(n_1, \ldots, n_s) \neq (0, \ldots, 0)$ such that
\[
H^0(X, {\mathcal O}_X(\sum_i n_iD_i)) \neq 0 \mbox{ \ \ and \ \ }
H^0(X, {\mathcal O}_X(-\sum_i n_iD_i)) \neq 0 .
\]
Then, there exists an effective Weil divisor $F_1$ that is linearly equivalent to $\sum_i n_iD_i$.
Since $\overline{D_1}$, \ldots, $\overline{D_s}$ are linearly independent over ${\mathbb Z}$ in the divisor class group ${\rm Cl}(X)$
by our assumption,
$F_1 \neq 0$.
In the same way, there exists a non-zero effective Weil divisor $F_2$ 
that is linearly equivalent to $-\sum_i n_iD_i$.
Then, the non-zero effective Weil divisor $F_1 + F_2$ is linearly equivalent to $0$.
This is a contradiction.
\qed

\begin{Remark}\label{rem}
\begin{rm}
Suppose that the assumptions (2) and (3) in Theorem~\ref{main} are satisfied. 
Set $R = R(X; D_1, \ldots, D_s)$.

Suppose that $D = a_1D_1 + \cdots + a_sD_s$ is an ample Cartier divisor, 
where $a_1, \ldots, a_s \in {\mathbb Z}$.
Set
\[
S = 
\bigoplus
_{n \geq 0} H^0(X, {\mathcal O}_X(nD))t^n \subset k(X)[t] .
\]
We have a ring homomorphism $S \rightarrow R$ defined by $ft^n \mapsto f(t_1^{a_1} \cdots t_s^{a_s})^n$
for $f \in H^0(X, {\mathcal O}_X(nD))$.
We think $S$ as a subring of $R$ by this ring homomorphism.
Here, $X = \proj(S)$.
%Let $S_+$ be the unique homogeneous maximal ideal of $S$.
Set $X' = X \setminus \Sing(X)
$ and $D'_i = D_i|_{X'}$ for $i = 1, \ldots, s$.

Replacing $D$ by $nD$ for a sufficiently large $n$ (if necessary),
we may assume that there  exist $f_1, \ldots, f_\ell \in H^0(X, {\mathcal O}_X(D))$
which satisfy the following two conditions:
\begin{itemize}
\item[(1)]
$X'  =  \bigcup
_{j = 1}^\ell D_+(f_jt)$.
\item[(2)]
$D'_i$ is a principal divisor on $D_+(f_jt)$ for any $i$ and $j$.
\end{itemize}

Let $(f_jt_1^{a_1} \cdots t_s^{a_s} \mid j ) R$ be the ideal of $R$
generated by $\{ f_jt_1^{a_1} \cdots t_s^{a_s} \mid j = 1, \ldots, \ell \}$.
Then, we can show that
\begin{equation}\label{ht}
\mbox{the height of  the ideal $(f_jt_1^{a_1} \cdots t_s^{a_s} \mid j ) R$ is bigger than 
or equal to two}
\end{equation}
as follows.
Let $H_1$ be the set of closed subvarieties of $X$ of codimension one.
For $V \in H_1$, we set
\begin{equation}\label{P_V}
P_V = 
\bigoplus_{(n_1, \ldots, n_s) \in {\mathbb Z}^s} H^0(X, {\mathcal O}_X(\sum_i n_iD_i - V))t_1^{n_1}\cdots t_s^{n_s}
\subset R ,
\end{equation}
that is, $P_V = M_{-V}$.
Then, as in p.~632 in \cite{EKW}, 
\[
\{ P_V \mid V \in H_1 \}
\]
coincides with the set of all the height one homogeneous prime ideals of $R$.
Here, recall that $P_V \cap S$ is equal to the defining ideal of $V$ in the ring $S$.
Since $V_+((f_jt \mid j )S)$ coincides with $\Sing(X)$,
the ideal $(f_jt_1^{a_1} \cdots t_s^{a_s} \mid j ) R$ is not contained in any
height one homogeneous prime ideal of $R$.
Here, recall that $\Sing(X)$ is a closed subset of $X$ of codimension 
bigger than or equal to $2$.

Thus, we know that the height of the ideal $(f_jt_1^{a_1} \cdots t_s^{a_s} \mid j ) R$ is bigger than 
or equal to two.
\end{rm}
\end{Remark}

\section{A proof of the main theorem using the twisted inverse functor}\label{proof}

We shall prove Theorem~\ref{main} in this section using the twisted inverse functor.
We shall prove
\[
\omega_{R(X; D_1, \ldots, D_s)}
\simeq M_{K_X}
\]
without assuming (1) in Theorem~\ref{main}.
If we remove the assumption~(1) in Theorem~\ref{main},
the ring $R(X; D_1, \ldots, D_s)$ may not be a local ${\mathbb Z}^s$-graded
$k$-domain.
In this case, we can not define $\omega_{R(X; D_1, \ldots, D_s)}$
using local cohomologies as in Definition~\ref{canonical}.
In this section, we shall give an alternative definition of  $\omega_{R(X; D_1, \ldots, D_s)}$
using the twisted inverse functor as in Definition~\ref{newdef}.
Of course, the both definitions coincide in the case where
$R(X; D_1, \ldots, D_s)$ is a local ${\mathbb Z}^s$-graded
$k$-domain as in Remark~\ref{coincides}.

\begin{Definition}\label{newdef}
Let $k$ be a field, $G$ a finite-type group scheme over $k$, and
$f:X\rightarrow \Spec k$ 
a $G$-scheme which is separated of finite type.
Then the dualizing complex $\Bbb I_X$ of $X$ is defined to be $f^!\O_{\Spec k}
$, 
where $f^!$ denotes 
the (equivariant) twisted inverse functor $f^!:D_{\Lqc}(G,\Spec k)
\rightarrow D_{\Lqc}(G,X)$, see  \cite[(20.5)]{LH} and \cite[Chapter~29]{LH}.
Assume that $X$ is non-empty and connected.
Then we define $s:=\inf\{i\mid H^i(\Bbb I_X)\neq 0\}$, and 
$\omega_X:=H^s(\Bbb I_X)$.
We call $\omega_X$ the ($G$-equivariant) canonical sheaf of $X$.
If%
%%% added 2010/11/15
, moreover, 
%%%
$X=\Spec A$ is affine, then $\Gamma(X,\omega_X)$ is denoted by $\omega_A$,
and is called the ($G$-equivariant) canonical module of $A$.
\end{Definition}

\begin{Remark}\label{canonical.rem}
Note that $\Bbb I_X\in D^b_{\Coh}(G,X)$, and $\omega_X\in \Coh(G,X)$.
Note also that if we forget the $G$-action, then $\Bbb I_X$ as an
object of $D^b(X)$ is the dualizing complex of $X$, and $\omega_X$ is
the canonical sheaf of $X$ in the usual sense.
If $X$ is Cohen--Macaulay, then $\omega_X[-s]\cong \Bbb I_X$ in 
$D(G,X)$, 
where $s=\inf\{i\mid H^i(\Bbb I_X)\neq 0\}$
(we use this fact freely later, for non-singular varieties).
If $X$ is equidimensional and $U$ is a $G$-stable open subset, then
$\omega_X|_U\cong\omega_U$.
In general, for a 
%%% added 2010/11/15
quasi-compact 
%%%
separated $G$-morphism $g:U\rightarrow X$,
the canonical map $u:\omega_X\rightarrow g_*g^*\omega_X$ is $(G,\O_X)$-linear.
If $X$ is normal, $g$ is the inclusion from a $G$-stable open subset $U$,
and the codimension of $X \setminus U$ in $X$ is at least two, then $u$ is an isomorphism
of coherent $(G,\O_X)$-modules.
So we have $\omega_X\cong g_*\omega_U$ this case.
\end{Remark}

\begin{Remark}\label{coincides-pre.rem}\rm
Let $k$ be a field, $H$ a finite-type $k$-group scheme, and
$G:=\Bbb G_m\times H$.
A $G$-algebra $R$ is nothing but a $\Bbb Z$-graded $k$-algebra 
$R=\bigoplus_{i\in\Bbb Z} R_i$ 
which is also an $H$-algebra such that each $R_i$ is an 
$H$-submodule of $R$, where the $\Bbb Z$-grading is given by the 
$\Bbb G_m$-action.
Let $R$ be a positively graded $G$-algebra.
That is, $R=\bigoplus_{i\in\Bbb Z}R_i$ is a $G$-algebra
such that $R_i=0$ for $i<0$ and $R_0=k$.
Assume that $R$ is a finitely generated domain.

Let $d=\dim R$ and $\Theta=\Spec R$.
Let $\theta$ be the unique $G$-stable closed point of $\Theta$,
corresponding to the unique graded maximal ideal $\frak m$ of $R$.
Then $\omega_\Theta=H^{-d}(\Bbb I_\Theta)$.
By the equivariant local duality \cite[(4.18)]{HO}, for any $\Bbb F
\in D_{\Coh}(G,\Theta)$, 
\[
\uH^i_\theta(\Bbb F)\cong \uHom_{\O_\Theta}(\uExt_{\O_\Theta}^{-i}
(\Bbb F,\Bbb I_\Theta),\E),
\]
where 
$\uHom_{\O_\Theta}$ denotes the sheaf-Hom in the category of
$(G,\O_\Theta)$-modules \cite[(2.24), Chapter~29]{LH},
and $\uExt^{-i}_{\O_\Theta}$ its right derived functor.
$\uH_\theta^i$ denotes the equivariant local cohomology 
\cite[(7.7)]{HO2}, and
$\E=\uH^0_\theta(\Bbb I_\Theta)$ is the $G$-sheaf of Matlis
of the $G$-local $G$-scheme $(\Theta,\theta)$ \cite{HO}.
Letting $\Bbb F=\O_\Theta$ and $i=d$,
\[
\uH_\theta^d(\O_\Theta)\cong\uHom_{\O_\Theta}(\omega_\Theta,\E).
\]
In other words,
\[
H^d_{\fm} (R)\cong \Hom_R(\omega_R,E),
\]
where $E=H^0(\Theta,\E)$.
By \cite[Lemma~5.4]{HO} and \cite[Remark~5.6]{HO}, 
it is easy to see that $(?)^\vee=
\Hom_R(?,E)$ is equivalent to ${}^*\Hom_k(?,k)$ 
as a functor from the full subcategory of the category of $(G,R)$-modules
consisting of modules whose degree $i$ component is a finite dimensional
$k$-vector space for each $i\in\Bbb Z$ (the $\Bbb Z$-grading is given
by the $\Bbb G_m$-action) to itself, 
and $(?)^\vee(?)^\vee$ is equivalent to the identity functor.
Thus $(H^d_{\fm} (R))^\vee\cong \omega_R$.
\end{Remark}

Let $H=\Bbb G_m^s$.
\iffalse
In general, for an $H$-scheme $X$, 
the category of quasi-coherent $(H,\O_X)$-modules $\Qch(H,X)$ is
identified with the category of $\Bbb Z^{s}$-graded $\O_X$-modules and
degree-preserving morphisms.
Assume that $X$ is Noetherian.
$\mathcal{M}\in\Qch(H,X)$ lies in $\Coh(H,X)$, the category of
coherent $(H,\O_X)$-modules, if and only if 
$\mathcal {M}$ is coherent as an $\O_X$-module.
\fi
%%% modified 2010/11/15
Let $\varphi:X\rightarrow Y$ be an affine $H$-morphism between $H$-schemes.
Assume that the action of $H$ on $Y$ is trivial.
Then $\varphi_*\O_X$ is a quasi-coherent $(H,\O_Y)$-algebra in the sense that
$\varphi_*\O_X$ is both 
a quasi-coherent $(H,\O_Y)$-module and an $\O_Y$-algebra, and
the unit map $\O_Y\rightarrow \varphi_*\O_X$ and the product
$\varphi_*\O_X\otimes_{\O_Y}\varphi_*\O_X\rightarrow \varphi_*\O_X$ 
are $(H,\O_Y)$-linear.
This is equivalent to say that $\varphi_*\O_X$ is a $\Bbb Z^s$-graded
$\O_Y$-algebra.
Conversely, for a given trivial $H$-scheme $Y$ and a quasi-coherent
$(H,\O_Y)$-algebra (that is, a $\Bbb Z^s$-graded $\O_Y$-algebra) $\A$, 
letting $X:=\Spec_Y \A$ and $\varphi:X\rightarrow Y$ the canonical map,
$\varphi:X\rightarrow Y$ is an affine $H$-morphism 
such that $\varphi_*\O_X\cong\A$.
In this case, a quasi-coherent $(H,\O_X)$-module $\M$ yields
a graded $\A$-module $\varphi_*\M$, and conversely, a graded $\A$-module
$\N$ determines a quasi-coherent $(H,\O_X)$-module $\M$ such that
$\varphi_*\M\cong\N$ uniquely, up to isomorphisms.
%%%

\begin{Remark}\label{coincides}
Let $H=\Bbb G_m^s$, and $G=\Bbb G_m\times H=\Bbb G_m^{s+1}$.

Let $R$ be a $\Bbb Z^s$-graded $k$-algebra.
Then $R$ is an $H$-algebra.
Assume moreover that $R$ is a local $\Bbb Z^s$-graded $k$-domain, see
Definition~\ref{local.def}.

Then the convex polyhedral cone in $\Bbb R^s$ generated by 
$\{\underline{a}\in\Bbb Z^s\mid R_{\underline{a}}\neq 0\}$ does not
contain a line.
So there is a linear function $\varphi:\Bbb Z^s\rightarrow \Bbb Z$ such
that $R=\bigoplus_{i\in\Bbb Z}R_i$ is a positively graded
$G$-algebra,
that is, $R_0=k$ and $R_i=0$ for $i<0$ and each $R_i$ is an $H$-submodule of
$R$, where 
$R_i=\bigoplus_{\varphi(\underline{a})=i}R_{\underline{a}}$.
Note that $R$ is assumed to be a finitely generated domain.

As in Remark~\ref{coincides-pre.rem}, $(H^d_{\fm} (R))^\vee\cong \omega_R$
as $(G,R)$-modules.
In particular, they are isomorphic as $(H,R)$-modules or 
$\Bbb Z^s$-graded $R$-modules.

This shows that our new definition of $\omega_R$ as a graded $R$-module is
consistent with Definition~\ref{canonical}.
\end{Remark}

We assume (2) and (3) in Theorem~\ref{main}.
Set $R = R(X; D_1, \ldots, D_s)$.
Take an ample Cartier divisor $D$, a subring $S$ and 
$f_1, \ldots, f_\ell \in H^0(X, {\mathcal O}_X(D))$ as in Remark~\ref{rem}.

We set
\[
Y = {\rm Spec}_{X'}\left(
\bigoplus_{(n_1, \ldots, n_s) \in {{\mathbb N}_0}^s} {\mathcal O}_{X'}(\sum_i n_iD'_i)
t_1^{n_1} \cdots t_s^{n_s}
\right) ,
\]
where ${\mathbb N}_0$ denotes the set of all non-negative integers.
Let $\pi : Y \rightarrow X'$ be the structure morphism.
Consider the open subscheme
\begin{equation}\label{Z}
Z = {\rm Spec}_{X'}\left(
\bigoplus_{(n_1, \ldots, n_s) \in {\mathbb Z}^s} {\mathcal O}_{X'}(\sum_i n_iD'_i)
t_1^{n_1} \cdots t_s^{n_s}
\right) 
\end{equation}
of $Y$.
Let $i : Z \rightarrow Y$ be the open immersion.

For each affine open subset $U$ of $X'$, we have a ring homomorphism
\[
R \longrightarrow
\bigoplus_{(n_1, \ldots, n_s) \in {\mathbb Z}^s}
H^0(U, {\mathcal O}_{X'}(\sum_in_iD'_i))t_1^{n_1} \cdots t_s^{n_s}
\]
induced by the restriction to the open set $U$ of $X'$.
Therefore, we have a natural morphism $j : Z \rightarrow \spec(R)$.

\begin{equation}\label{schemes}
\begin{array}{ccccc}
& & Z & \stackrel{i}{\longrightarrow} & Y \\
& & \phantom{\scriptstyle j}\downarrow{\scriptstyle j} &  & 
\phantom{\scriptstyle \pi}\downarrow{\scriptstyle \pi}  \\
\spec(R) \setminus V((f_jt_1^{a_1} \cdots t_s^{a_s} \mid j ) R) & \subset & \spec(R) & & X'
\end{array}
\end{equation}

\begin{Lemma}\label{lem1}
The morphism $j : Z \rightarrow \spec(R)$ coincides with the open immersion
$\spec(R) \setminus V((f_jt_1^{a_1} \cdots t_s^{a_s} \mid j ) R)  \subset  \spec(R)$.
\end{Lemma}

\proof
There exist $\alpha_{ij}$'s in $k(X)^\times$ such that
\[
R[(f_jt_1^{a_1} \cdots t_s^{a_s})^{-1}] 
= \Gamma(D_+(f_jt), {\mathcal O}_X)[(\alpha_{1j}t_1)^{\pm 1}, \ldots, (\alpha_{sj}t_s)^{\pm 1}] 
\]
for each $j$.
The ring corresponding to the affine open set $(\pi i)^{-1}(D_+(f_jt))$ just coincides with the above ring.
Thus we obtain 
\[
Z = \spec(R) \setminus V((f_jt_1^{a_1} \cdots t_s^{a_s} \mid j ) R) .
\]
\qed

The group $({\mathbb G}_m)^s$ acts on all the schemes in diagram~(\ref{schemes}).
The action on $X'$ is trivial.
All the morphisms in diagram~(\ref{schemes}) are compatible with the group action.

As in Definition~\ref{newdef}, we can define the canonical sheaves with group action for the schemes in 
diagram~(\ref{schemes}).
That is, the canonical sheaves for all the schemes in diagram~(\ref{schemes})
are ${\mathbb Z}^s$-graded.

Then, we have $\omega_R|_Z = \omega_Z$ by the compatibility with open immersions
(see Remark~\ref{canonical.rem}).
Since the height of $(f_jt_1^{a_1} \cdots t_s^{a_s} \mid j ) R$ is bigger than one as in Remark~\ref{rem},
we have an isomorphism
\[
H^0(Z, \omega_Z) = \omega_R
\]
which is compatible with the group action (see Remark~\ref{canonical.rem}).

On the other hand, $\pi : Y \rightarrow X'$ is a smooth morphism of relative dimension $s$ that
is compatible with the group action.
Then, the sheaf of differentials $\Omega_{Y/X'}$ naturally has a group action, that is,
it has a structure of a ${\Bbb Z}^s$-graded sheaf.
We have 
\[
\bigwedge^s \Omega_{Y/X'} \simeq \pi^*{\mathcal O}_{X'}(\sum_iD'_i)
(-1, \ldots, -1) ,
\]
where $(-1, \ldots, -1)$ denotes the shift of degree.

Then, by Theorem~28.11 in \cite{LH}, we have
\begin{eqnarray*}
\omega_Y & = & \bigwedge^s \Omega_{Y/X'} \otimes_{{\mathcal O}_Y} \pi^*\omega_{X'} \\
& \simeq & \pi^*{\mathcal O}_{X'}(\sum_iD'_i)(-1, \ldots, -1) \otimes_{{\mathcal O}_Y} 
\pi^*{\mathcal O}_{X'}(K_{X'}) \\
& = & \pi^*{\mathcal O}_{X'}(\sum_iD'_i + K_{X'})(-1, \ldots, -1) .
\end{eqnarray*}

So, we have
\[
\omega_Z = i^* \omega_Y = (\pi i)^*{\mathcal O}_{X'}(\sum_iD'_i + K_{X'})
(-1, \ldots, -1).
\]

Then,
\begin{eqnarray*}
\omega_R & = & H^0(Z, \omega_Z) \\
& = & H^0(X', (\pi i)_*\omega_Z) \\
& = & H^0(X', (\pi i)_*(\pi i)^*{\mathcal O}_{X'}(\sum_iD'_i + K_{X'})
(-1, \ldots, -1)) .
\end{eqnarray*}

By the equivariant projection formula 
(Lemma~26.4 in \cite{LH}),
we have
\begin{eqnarray*}
& & (\pi i)_*(\pi i)^*{\mathcal O}_{X'}(\sum_iD'_i + K_{X'})
(-1, \ldots, -1) \\
& \simeq & \left(
{\mathcal O}_{X'}(\sum_iD'_i + K_{X'}) \otimes_{{\mathcal O}_{X'}} (\pi i)_*{\mathcal O}_Z
\right) (-1, \ldots, -1) \\
& = & \left(
{\mathcal O}_{X'}(\sum_iD'_i + K_{X'}) \otimes_{{\mathcal O}_{X'}} 
\left[
\bigoplus_{(n_1, \ldots, n_s) \in {\mathbb Z}^s} {\mathcal O}_{X'}(\sum_i n_iD'_i)
\right]
\right) (-1, \ldots, -1) \\
& = & \left(
\bigoplus_{(n_1, \ldots, n_s) \in {\mathbb Z}^s} {\mathcal O}_{X'}(\sum_i (n_i+1)D'_i + K_{X'})
\right) (-1, \ldots, -1) \\
& = & 
\bigoplus_{(n_1, \ldots, n_s) \in {\mathbb Z}^s} {\mathcal O}_{X'}(\sum_i n_iD'_i + K_{X'}) .
\end{eqnarray*}

Thus, we obtain
\begin{eqnarray*}
\omega_R & = & 
\bigoplus_{(n_1, \ldots, n_s) \in {\mathbb Z}^s} H^0(X', {\mathcal O}_{X'}(\sum_i n_iD'_i + K_{X'})) \\
& = &
\bigoplus_{(n_1, \ldots, n_s) \in {\mathbb Z}^s} H^0(X, {\mathcal O}_{X}(\sum_i n_iD_i + K_{X})) .
\end{eqnarray*}

We have completed the proof of Theorem~\ref{main}.

\section{Another proof of Theorem~\ref{main}}\label{another}

In this section, we give another proof to Theorem~\ref{main} without using the twisted inverse functor.
In this section, we have to assume that the scheme $X'$ in Remark~\ref{rem} is smooth over $k$.
Note that it is automatically satisfied if the base field $k$ is perfect.

The idea of this proof is based on Lemma~13 in \cite{P}.
In the proof in this section, we have to assume that Theorem~\ref{main} is true if $s = 1$.

\begin{Remark}\label{Serre}
\begin{rm}
In the case where $s = 1$, Theorem~\ref{main} is proved using Serre duality.
Assume that $s = 1$ and $a_1D_1$ is a very ample Cartier divisor for some $a_1 > 0$.
Put $R = R(X;D_1)$.
Then, we have
\begin{eqnarray*}
& & \bigoplus_{n_1 \in {\mathbb Z}} H^0(X, {\mathcal O}_{X}(n_1D_1 + K_X))  = 
\bigoplus_{n_1 \in {\mathbb Z}} 
{\rm Hom}_{{\mathcal O}_{X}}({\mathcal O}_{X}(-n_1D_1), \omega_X) \\
& = & 
\bigoplus_{n_1 \in {\mathbb Z}} 
{\rm Hom}_{k}(H^d(X, {\mathcal O}_{X}(-n_1D_1)), k) 
= 
\bigoplus_{n_1 \in {\mathbb Z}} 
{\rm Hom}_{k}(H^{d+1}_{R_+}(R)_{-n_1}, k) \\
& = & {}^*{\rm Hom}_{k}(H^{d+1}_{R_+}(R), k) = \omega_R .
\end{eqnarray*}
\end{rm}
\end{Remark}

First, we shall prove a basic fact on class groups as follows.
The authors guess that it is well-known, however we do not know an adequate reference.
Therefore, we shall give a proof here.

\begin{Lemma}\label{lemma4.1}
Let $U$ be a normal scheme.
Let $F_1$, \ldots, $F_t$ be Cartier divisors on $U$.
Put
\[
W = \spec_U
\left(
\bigoplus_{m_1, \ldots, m_t \in {\mathbb Z}} {\mathcal O}_U(\sum_{j=1}^t m_jF_j)
\right) .
\]
Let $\pi : W \rightarrow U$ be the structure morphism.

Then, the sequence
\[
0 \longrightarrow {\mathbb Z}\overline{F_1} + \cdots + {\mathbb Z}\overline{F_t}
\longrightarrow {\rm Cl}(U) \stackrel{\pi^*}{\longrightarrow} {\rm Cl}(W) 
\longrightarrow 0
\]
is exact.
\end{Lemma}

\proof
Recall that, since $\pi : W \rightarrow U$ is a flat morphism,
the pull-back map 
\[
\pi^* : {\rm Cl}(U) \longrightarrow {\rm Cl}(W)
\]
is defined as in \cite{F}.

Put
\begin{eqnarray*}
E & = & \spec_U
\left(
\bigoplus_{m_1, \ldots, m_t \geq 0} {\mathcal O}_U(\sum_{j=1}^t m_iF_i) 
\right)
= \spec_U
\left(
{\rm Sym}_{{\mathcal O}_U}\left( 
{\mathcal O}_U(F_1) \oplus \cdots \oplus {\mathcal O}_U(F_t)
\right)
\right)
\\
E_j & = & \spec_U
\left(
\bigoplus_{m_j \geq 0} {\mathcal O}_U(m_jF_j) 
\right)
= \spec_U
\left(
{\rm Sym}_{{\mathcal O}_U}\left( 
{\mathcal O}_U(F_j)
\right)
\right)
\end{eqnarray*}
for $j = 1, \ldots, t$.
Let 
\begin{eqnarray*}
p & : & E \longrightarrow U 
\\
p_j & : & E_j \longrightarrow U
\\
q_j & : & E \longrightarrow E_j
\end{eqnarray*}
be the natural morphisms.
Remark that $p = p_jq_j$ for any $j$.
By Theorem~3.3 (a) in \cite{F}, the flat pull-back map
\[
p^* : {\rm Cl}(U) \longrightarrow {\rm Cl}(E)
\]
is an isomorphism.
Let 
\[
s_j : U \longrightarrow E_j
\]
be the $0$-section of the vector bundle $p_j : E_j \rightarrow U$.
Then, 
\[
W = E \setminus \bigcup_j q_j^{-1}(s_j(U)) .
\]
By Proposition~1.8 in \cite{F}, we have the following exact sequence:
\[
\begin{array}{ccccccccc}
0 & \longrightarrow & \sum_j {\mathbb Z}\overline{q_j^{-1}(s_j(U))}
& \longrightarrow & {\rm Cl}(E) & \longrightarrow & {\rm Cl}(W) & \longrightarrow & 0 \\
& & & & \hphantom{\scriptstyle p^*} \uparrow \scriptstyle p^* & & & & \\
& & & & {\rm Cl}(U) & & & & 
\end{array}
\]
Then $s_j(U)$ is linearly equivalent to $p_j^*(-F_j)$ by Example~3.3.2 in \cite{F}.
Therefore, $q_j^{-1}(s_j(U))$ is linearly equivalent to $-p^*(F_j)$.
Thus, we have the desired exact sequence.
\qed

\begin{Lemma}\label{lem2}
Suppose that (2) in Theorem~\ref{main} 
%%% modified 2011/03/10
is
%%% are 
satisfied.

Then, we have the following:
\begin{enumerate}
\item
The set
\[
\{ M_F \mid \mbox{$F$ is a Weil divisor on $X$} \}
\]
just coincides with the set of divisorial fractional ideals
which are ${\mathbb Z}^s$-graded $R(X; D_1, \ldots, D_s)$-submodules of 
$k(X)[t_1^{\pm 1}, \ldots, t_s^{\pm 1}] $.
\item
For Weil divisors $F_1$ and $F_2$,
$F_1$ is linearly equivalent to $F_2$ if and only if 
$M_{F_1}$ is isomorphic to $M_{F_2}$ as a ${\mathbb Z}^s$-graded module.
\item
Further, assume that $R(X; D_1, \ldots, D_s)$ is Noetherian.
Any non-zero finitely generated ${\mathbb Z}^s$-graded
reflexive $R(X; D_1, \ldots, D_s)$-module of rank one is isomorphic to 
$M_G$ as a ${\mathbb Z}^s$-graded module for some Weil divisor $G$ on $X$.
\end{enumerate}
\end{Lemma}

\proof
Set $R = R(X; D_1, \ldots, D_s)$.

First, we prove (1).
Let $F$ be a Weil divisor on $X$.
We shall prove that $M_F$ is a divisorial fractional ideal.

Since there exists an ample Cartier divisor in ${\mathbb Z}D_1+ \cdots +{\mathbb Z}D_s$,
we can find $a t_1^{m_1} \cdots t_s^{m_s}$ such that ${\rm div}_X(a) + \sum_im_iD_i - F$
is an effective divisor, 
where ${\rm div}_X(a)$ is the principal Weil divisor corresponding to $a \in k(X)^\times$.
By definition, it is easy to check that, for $a \in k(X)^\times$,
\begin{equation}\label{siki}
(a t_1^{m_1} \cdots t_s^{m_s} ) M_F = M_{F - {\rm div}_X(a) - \sum_im_iD_i} .
\end{equation}
We have only to show that $M_{F - {\rm div}_X(a) - \sum_im_iD_i}$
is a divisorial fractional ideal.

We define $H_1$ and $P_V$ as in Remark~\ref{rem}.
Set 
\[
F - {\rm div}_X(a) - \sum_im_iD_i = - \sum_{j = 1}^u\ell_jV_j, 
\]
where $V_1$, \ldots, $V_u$ are distinct 
elements in $H_1$.
By the argument as above, all of $\ell_j$'s are positive integers.

For $V \in H_1$,  we define $R_V$ as in 
%%% modified 2011/03/11
%629p 
p.~629
in \cite{EKW}.
(Here, we use the symbol "$V$" instead of  "$F$" in \cite{EKW}.)
In this case, we obtain
\[
R_{P_V} = (R_V)_{\alpha_V R_V}
\]
as in p.~632 in \cite{EKW}.
Therefore, for any $\ell > 0$,
we have
\[
{P_V}^{(\ell)} = {P_V}^\ell R_{P_V} \cap R
= {\alpha_V}^\ell (R_V)_{\alpha_V R_V} \cap R
= M_{-\ell V} .
\]
Then,
\[
M_{- \sum_j \ell_jV_j} = \bigcap_j M_{- \ell_jV_j} =  \bigcap_j {P_{V_j}}^{(\ell_j)} .
\]
Since $P_{V_j}$'s are homogeneous prime ideals of $R$ of height one,
$M_{- \sum_j\ell_jV_j}$ is a divisorial fractional ideal.

Conversely, let $N$ be a divisorial fractional ideal that is a ${\mathbb Z}^s$-graded
$R$-submodule of $k(X)[t_1^{\pm 1}, \ldots, t_s^{\pm 1}] $.
Using (\ref{siki}), it is sufficient to prove that $N$ coincides with $M_F$
for some Weil divisor $F$ in the case where $N \subset R$.
Then, $N$ coincides with an intersection of symbolic powers of homogeneous prime ideals of height one.
Therefore, there exist $V_1, \ldots, V_u \in H_1$ and positive integers $\ell_1, \ldots, \ell_u$
such that
\[
N = \bigcap_j {P_{V_j}}^{(\ell_j)} = \bigcap_j M_{- \ell_jV_j} = M_{- \sum_j\ell_jV_j} .
\]

(2) and (3) are easily verified, so we omit the proofs.
\qed

\begin{Remark}\label{divisorial}
\begin{rm}
Assume (2) in Theorem~\ref{main}.
Put $R = R(X; D_1, \ldots, D_s)$.
The map $p$ in (\ref{ekw}) is constructed as follows.

Let ${\rm HDiv}(R)$ be the free abelian group generated 
by all the homogeneous prime divisors of $\spec(R)$, i.e.,
\[
{\rm HDiv}(R)
= 
\bigoplus_{V \in H_1(X)} {\Bbb Z} [\spec(R/P_V)] ,
\]
where $ [\spec(R/P_V)]$ denotes the generator 
corresponding to the closed subscheme $\spec(R/P_V)$.
Here, recall that 
\[
\{ P_V \mid V \in H_1(X) \}
\]
coincides with the set of all the homogeneous prime ideals of $R$
of height one as in p.~632 in \cite{EKW}.

We define 
\[
\xi' : {\rm Div}(X) \rightarrow {\rm HDiv}(R)
\]
by
\[
\xi' \left( \sum_{V \in H_1(X)} n_VV \right) = \sum_{V \in H_1(X)} n_V [\spec(R/P_V)] .
\]
(This map $\xi'$ is equal to $\xi$ in p631 in \cite{EKW}
if we identify $[\spec(R/P_V)]$ with $P_V$.)
As in \cite{EKW}, $\xi'$ induces the map
\[
\varphi' : {\rm Cl}(X) \rightarrow {\rm A}^1(\spec(R))
\]
satisfying
\[
\varphi' \left( \overline{\sum_{V \in H_1(X)} n_VV} \right) = 
\overline{\sum_{V \in H_1(X)} n_V [\spec(R/P_V)]} ,
\]
where ${\rm A}^1(\spec(R))$ denotes the Chow group of $\spec(R)$ of codimension one.
Then, we have the following exact sequence
\[
0 \longrightarrow {\mathbb Z}\overline{D_1} + \cdots + {\mathbb Z}\overline{D_s}
\longrightarrow {\rm Cl}(X)
\stackrel{\varphi'}{\longrightarrow} {\rm A}^1(\spec(R)) \longrightarrow 0 .
\]

Let ${\rm Cl}(R)$ be the set of isomorphism classes of divisorial fractional ideals
of $R$.
Then we have the isomorphism
\[
i : {\rm A}^1(\spec(R)) \longrightarrow {\rm Cl}(R)
\]
defined by
\[
i \left( \overline{\sum_{V \in H_1(X)} n_V [\spec(R/P_V)]} \right) = 
-\sum_{V \in H_1(X)} n_V \overline{P_V} =
\overline{ M_{\sum_{V \in H_1(X)} n_VV} } ,
\]
where $\overline{P_V}$ (resp.\ $\overline{ M_{\sum_{V \in H_1(X)} n_VV} }$) denotes
the isomorphism class containing $P_V$ (resp.\ $M_{\sum_{V \in H_1(X)} n_VV}$).

We set $p = i \varphi'$.
Then,
\[
p : {\rm Cl}(X) \rightarrow {\rm Cl}(R)
\]
is the map which satisfy $p(\overline{F}) = \overline{ M_{F}}$ for each Weil divisor $F$ on $X$,
and we have the exact sequence (\ref{ekw}).
\end{rm}
\end{Remark}

We now start to prove Theorem~\ref{main} in the case where $s \geq 2$.
Suppose that all the assumptions in Theorem~\ref{main} are satisfied.
We assume that $X'$ in Remark~\ref{rem} is smooth over $k$.

Since $\omega_R$ is a finitely generated ${\mathbb Z}^s$-graded
reflexive $R$-module of rank one, there exists a Weil divisor $G$
such that $\omega_R$ is isomorphic to $M_G$ as a ${\mathbb Z}^s$-graded module
as in Lemma~\ref{lem2}.

We want to show that $G$ is linearly equivalent to $K_X$.
Assume the contrary, that is, $\overline{G-K_X} \neq 0$ in ${\rm Cl}(X)$.

Then, we can choose $F_1, \ldots, F_s \in {\mathbb Z}D_1+ \cdots +{\mathbb Z}D_s$
satisfying the following three conditions:
\begin{enumerate}
\item
${\mathbb Z}F_1 + \cdots + {\mathbb Z}F_s = {\mathbb Z}D_1+ \cdots +{\mathbb Z}D_s$.
\item
If $b_1F_1 + \cdots + b_sF_s$ is linearly equivalent to a non-zero effective divisor,
then $b_s > 0$.
\item
$\overline{G-K_X} \not\in {\mathbb Z}\overline{F_1} + \cdots + {\mathbb Z}\overline{F_{s-1}}$.
\end{enumerate}
Here, recall that $s$ is bigger than $1$, and
$F_1$, \ldots, $F_s$ are linearly independent over ${\mathbb Z}$
by the assumption (1) in Theorem~\ref{main}.

We define a map
\[
\varphi : {\mathbb Z}D_1+ \cdots +{\mathbb Z}D_s \longrightarrow {\mathbb Z}
\]
by 
\[
\varphi(b_1F_1 + \cdots + b_sF_s) = b_s .
\]
The kernel of $\varphi$ is equal to $ {\mathbb Z}F_1+ \cdots +{\mathbb Z}F_{s-1}$.

We think that $R$ is a ${\mathbb Z}$-graded ring by 
\[
{\rm deg}(a t_1^{n_1} \cdots t_s^{n_s}) = 
\varphi(n_1D_1 + \cdots + n_sD_s) .
\]
Set $T = \proj(R)$.
Take an ample Cartier divisor $D = \sum_ia_iD_i$, 
a subring $S$ and $f_1, \ldots, f_\ell$ as in Remark~\ref{rem}.
Set $Q = T \setminus V_+((f_jt_1^{a_1} \cdots t_s^{a_s} \mid j ) R)$.
In this section, for an 
$\ell$-dimensional normal algebraic variety $W$ over $k$ which is smooth over $k$
in codimension one, we define
\[
\omega_W = \left(  \bigwedge^\ell \Omega_{W/k} \right)^{**} .
\]
However, for $R$, we define $\omega_R$ as in Definition~\ref{canonical}.

Then, it is easy to see that $Q$ coincides with 
\[
\spec_{X'}\left(
\bigoplus_{m_1, \ldots, m_{s-1} \in {\mathbb Z}}
{\mathcal O}_{X'}(\sum_{j=1}^{s-1}m_jF'_j)
\right)   ,
\]
where $F'_j = F_j|_{X'}$.
Since $X \setminus X' = 
%%% modified by Hashimoto 2010/06/23
%{\rm Sing}
\Sing
(X)$,
the natural restriction
\begin{equation}\label{con}
{\rm Cl}(X) \longrightarrow {\rm Cl}(X')
\end{equation}
is an isomorphism.
Let $\pi : Q \rightarrow X'$ be the structure morphism.  

Then, by Lemma~\ref{lemma4.1}, we have the following exact sequence:
\begin{equation}\label{ex}
0 \longrightarrow {\mathbb Z}\overline{F'_1} + \cdots + {\mathbb Z}\overline{F'_{s-1}}
\longrightarrow {\rm Cl}(X') \stackrel{\pi^*}{\longrightarrow} {\rm Cl}(Q) \longrightarrow 0
\end{equation}

Since $\dim Q = d + s -1$, we have
\[
\omega_T|_Q = \omega_Q = \bigwedge^{d+s-1} \Omega_{Q/k} .
\]
Since $X'$ is smooth over $k$ and $\pi : Q \rightarrow X'$ is a smooth morphism,
the sequence
\[
0 \longrightarrow \pi^*\Omega_{X'/k} \longrightarrow \Omega_{Q/k}
\longrightarrow \Omega_{Q/X'} \longrightarrow 0
\]
is exact.
Therefore we have
\begin{eqnarray}
\omega_Q & = & \bigwedge^{d+s-1} \Omega_{Q/k} \nonumber \\
& = & \bigwedge^d\pi^*\Omega_{X'/k} \otimes_{{\mathcal O}_Q}
\bigwedge^{s-1} \Omega_{Q/X'} \nonumber \\
& = & \pi^* \bigwedge^d\Omega_{X'/k} \otimes_{{\mathcal O}_Q}
\pi^* {\mathcal O}_{X'}(F'_1 + \cdots + F'_{s-1}) \nonumber \\
& = & \pi^* \omega_{X'} \otimes_{{\mathcal O}_Q} {\mathcal O}_Q \nonumber \\
& = & \pi^* {\mathcal O}_{X'}(K_{X'})  \label{K_X}
\end{eqnarray}
since $\pi^* {\mathcal O}_{X'}(F'_1 + \cdots + F'_{s-1}) \simeq {\mathcal O}_Q$ 
by the exact sequence~(\ref{ex}).

On the other hand, we have
\[
\omega_T = \left(
\bigoplus_{(n_1, \ldots, n_s) \in {\mathbb Z}^s}
H^0(X, {\mathcal O}_X(\sum_in_iD_i + G))
\right) ^\sim ,
\]
where the right hand side is the coherent ${\mathcal O}_T$-module
associated with the ${\mathbb Z}$-graded module 
\[
\bigoplus_{(n_1, \ldots, n_s) \in {\mathbb Z}^s}
H^0(X, {\mathcal O}_X(\sum_in_iD_i + G)) ,
\]
that is the canonical module of the ${\mathbb Z}$-graded ring $R$.
(Remark that we used the fact that Theorem~\ref{main} is true if $s = 1$.)
Therefore, we obtain
\[
\omega_Q = \omega_T|_Q  = \left(
\bigoplus_{m_1, \ldots, m_{s-1} \in {\mathbb Z}}
{\mathcal O}_{X'}(\sum_{j=1}^{s-1}m_jF'_j + G')
\right) ^\sim ,
\]
where $G' = G|_{X'}$, and the right hand side is a coherent ${\mathcal O}_Q$-module
associated with 
\[
\bigoplus_{m_1, \ldots, m_{s-1} \in {\mathbb Z}}
{\mathcal O}_{X'}(\sum_{j=1}^{s-1}m_jF'_j + G') ,
\]
that is a sheaf of modules over the sheaf of algebras
\[
\bigoplus_{m_1, \ldots, m_{s-1} \in {\mathbb Z}}
{\mathcal O}_{X'}(\sum_{j=1}^{s-1}m_jF'_j)
\]
on $X'$.
Therefore, we have
\begin{equation}\label{G}
\omega_Q = \pi^* {\mathcal O}_{X'}(G') .
\end{equation}

Here, for Weil divisors $E_1$ and $E_2$ on $X'$, we know
\[
\pi^* {\mathcal O}_{X'}(E_1) \simeq \pi^* {\mathcal O}_{X'}(E_2) 
\Longleftrightarrow 
\overline{E_1} \equiv \overline{E_2} 
\mod {\mathbb Z}\overline{F'_1} + \cdots + {\mathbb Z}\overline{F'_{s-1}} 
\]
by the exact sequence~(\ref{ex}).

Thus, by (\ref{K_X}) and (\ref{G}),
we obtain 
\[
\overline{G'-K_{X'}} \in {\mathbb Z}\overline{F'_1} + \cdots + {\mathbb Z}\overline{F'_{s-1}}
.
\]
Since the restriction (\ref{con}) is an isomorphism, we know that 
\[
\overline{G-K_{X}} \in {\mathbb Z}\overline{F_1} + \cdots + {\mathbb Z}\overline{F_{s-1}}
\]
in ${\rm Cl}(X)$.
It is a contradiction.

We have completed the proof of Theorem~\ref{main}.

\begin{Remark}\label{4.3}
\begin{rm}
Suppose that all the assumptions in Theorem~\ref{main} are satisfied.
We choose $D$, $S$ as in Remark~\ref{rem}.
Put $d = \dim X$. 
Then, $\dim R = d + s$ as in \cite{EKW}.
Let ${\mathfrak m}$ be the unique maximal homogeneous ideal of $R$.
Let $S_+$ be the maximal ideal of $S$ generated by all the homogeneous elements
of positive degree.
Then, we have,
\begin{eqnarray*}
H^{d+s}_{{\mathfrak m}}(R) & = & {}^*{\rm Hom}_k(\omega_R,k) \\
& = & \bigoplus_{(n_1, \ldots, n_s) \in {\mathbb Z}^s}
{\rm Hom}_k(H^0(X, {\mathcal O}_X(-\sum_in_iD_i + K_X)), k) \\
& = & \bigoplus_{(n_1, \ldots, n_s) \in {\mathbb Z}^s}
H^d(X, {\mathcal O}_X(\sum_in_iD_i)) \\
& = & 
H^{d+1}_{S_+}(R) 
\end{eqnarray*}
by Serre duality.

Remark that there are many examples such that
\[
H^{d+s-1}_{{\mathfrak m}}(R) \neq H^{d}_{S_+}(R) .
\]
\end{rm}
\end{Remark}

\section{Some examples}\label{sectionex}

In this section, we give some examples.

\begin{Example}\label{ex1}
\begin{rm}
Let $B = k[x,y,z]$ be a weighted polynomial ring over a field $k$
with ${\rm deg}(x) = a$, ${\rm deg}(y) = b$ and ${\rm deg}(z) = c$, 
where $a$, $b$, $c$ are pairwise coprime positive integers.

Let $P$ be the kernel of the $k$-algebra homomorphism
$k[x,y,z] \longrightarrow k[t]$ defined by $x \mapsto t^a$, $y \mapsto t^b$, $z \mapsto t^c$.

Let $\pi : X \rightarrow {\rm Proj}(B)$ be the blow-up at $V_+(P)$.
Let $A$ be an integral Weil divisor on $X$ satisfying 
${\mathcal O}_X(A) = \pi^*{\mathcal O}_{{\rm Proj}(B)}(1)$. 
Put $E = \pi^{-1}(V_+(P))$.
In this case, ${\rm Cl}(X)$ is freely generated by $\overline{A}$ and $\overline{E}$.
We have
\[
K_X = E - (a+b+c)A .
\]

Then, $R(X; -E, A)$ coincides with the extended symbolic Rees ring
\[
R = k[x,y,z,t^{-1}, Pt, P^{(2)}t^2, P^{(3)}t^3, \ldots ] \subset k[x,y,z,t,t^{-1}] .
\]
$R$ is a ${\mathbb Z}^2$-graded ring with ${\rm deg}(x) = (0,a)$, 
${\rm deg}(y) = (0,b)$, ${\rm deg}(z) = (0,c)$, ${\rm deg}(t) = (1,0)$.
We know
\[
\omega_R = R(-1, -(a+b+c))  .
\]

For positive integers $\alpha$ and $\beta$,
we define
\[
R^{(\alpha, \beta)} = 
\bigoplus_{m_1, m_2 \in {\mathbb Z}} R_{(\alpha m_1, \beta m_2)} .
\]
Here, remark that $R^{(\alpha, \beta)} = R(X; -\alpha E, \beta A)$.
Therefore, $\omega_{R^{(\alpha, \beta)}}$ is an $R^{(\alpha, \beta)}$-free module
if and only if $\alpha = 1$ and $\beta | (a+b+c)$.
\end{rm}
\end{Example}

\begin{Example}\label{ex2}
\begin{rm}
Let $k$ be a field and $t_1$, \ldots, $t_r$ be rational distinct closed points
in ${\mathbb P}^n_k$.
Let $I_i$ be the defining ideal of $t_i$ in the homogeneous coordinate ring
$B = k[x_0, \ldots, x_n]$.
Recall that $I_i$ is generated by linearly independent $n$ linear forms.
Therefore, for any $\ell > 0$,
\[
I_i^\ell = (I_i^\ell)^{\rm sat} = I_i^{(\ell)},
\]
where $(I_i^\ell)^{\rm sat}$ is the {\em saturation} of the ideal $I_i^\ell$,
and $I_i^{(\ell)}$ is the $\ell$-th symbolic power of $I_i$.

Let $m_1$, \ldots, $m_r$ be positive integers.
Put
\[
I = I_1^{m_1} \cap \cdots \cap I_r^{m_r} .
\]

Let $\pi : X \rightarrow {\mathbb P}^n_k$ be the blow-up at $t_1$, \ldots, $t_r$.
Put $E_i = \pi^{-1}(t_i)$ for $i = 1, \ldots, \ell$.
Let $A$ be a Weil divisor on $X$ that satisfies ${\mathcal O}_X(A) = 
\pi^*{\mathcal O}_{{\mathbb P}^n_k}(1)$.
Then, ${\rm Cl}(X)$ is freely generated by $\overline{A}$,  $\overline{E_1}$,  
\ldots, $\overline{E_r}$.
In this case, we have
\[
K_X = (n-1)(E_1 + \cdots + E_r) - (n+1)A .
\]

Put
\[
R = R(X; -m_1E_1 - \cdots - m_rE_r, A) .
\]
Here, $-m_1E_1 - \cdots - m_rE_r + m A$ is an ample Cartier divisor for $m \gg 0$.
Then, $R$ coincides with the extended symbolic Rees ring of $I$, that is,
\[
R = B[t^{-1}, It, I^{(2)}t^2, I^{(3)}t^3, \ldots] \subset B[t, t^{-1}] ,
\]
where
\[
I^{(\ell)} = I_1^{\ell m_1} \cap \cdots \cap I_r^{\ell m_r} .
\]

Assume that $R$ is Noetherian.
Then, by Corollary~\ref{k}, we know 
\begin{eqnarray*}
\omega_R \simeq R & \Longleftrightarrow & 
\overline{K_X} \in {\mathbb Z}\overline{(-m_1E_1 - \cdots - m_rE_r)} + {\mathbb Z}\overline{A} \\
& \Longleftrightarrow & 
\left\{
\begin{array}{l}
m_1 = \cdots = m_r \\
m_1 \mid (n-1) .
\end{array}
\right.
\end{eqnarray*}
\end{rm}
\end{Example}

\begin{Example}\label{ex3}
\begin{rm}
Here, we give an example that the equality~(\ref{seigen}) is not satisfied
if we remove the assumption
that $L_2$ contains an ample Cartier divisor.

Set $X = {\mathbb P}^1_k \times {\mathbb P}^1_k$,
$A_1 = {\mathbb P}^1_k \times (1:0)$ and
$A_2 = (1:0) \times {\mathbb P}^1_k$.
Then, ${\rm Cl}(X)$ is freely generated by $\overline{A_1}$ and $\overline{A_2}$.
We know
\[
{\rm Cox}(X) = R(X;A_1,A_2)=k[x_0,x_1,y_0,y_1] ,
\]
where ${\rm deg}(x_0) = {\rm deg}(x_1) = (1,0)$ and
${\rm deg}(y_0) = {\rm deg}(y_1) = (0,1)$.

Then, we have
\[
\omega_{R(X;A_1,A_2)} = R(X;A_1,A_2)(-2,-2)
\]
since $K_X = -2A_1-2A_2$.

For $a, b \in {\mathbb Z}$, we define
\[
S_{a, b} = R(X; aA_1 + bA_2) 
\]
and
\[
L_{a, b} = {\mathbb Z}(aA_1 + bA_2) \subset {\mathbb Z}A_1 + {\mathbb Z}A_2 .
\]
Then, we have 
\[
S_{a, b} = R(X; A_1, A_2) |_{L_{a, b}} .
\]
In this case, $aA_1 + bA_2$ is ample if and only if $a > 0$ and $b > 0$.

Therefore, if $a > 0$ and $b > 0$, then 
\[
\omega_{S_{a, b}} = \omega_{R(X; A_1, A_2)}|_{L_{a, b}}
\]
by (\ref{seigen}).
However, we have
\[
\omega_{R(X; A_1, A_2)}|_{L_{1, 0}} = 0 \neq \omega_{S_{1, 0}} .
\]

In the case where $a > 0$ and $b > 0$,
$S_{a, b}$ is the Segre product of $k[x_0,x_1]^{(a)}$
and $k[y_0,y_1]^{(b)}$.
Here, $k[x_0,x_1]^{(a)}$ is the $a$-th Veronese subring of $k[x_0,x_1]$,
that is,
\[
k[x_0,x_1]^{(a)} = \bigoplus_{n \geq 0} k[x_0,x_1]_{na} .
\]
Then, we have
\[
\omega_{S_{a, b}} = \bigoplus_{n > 0} \left( k[x_0,x_1]_{na-2} \otimes_k k[y_0,y_1]_{nb-2} \right)
= \omega_{k[x_0,x_1]^{(a)}} \# \omega_{k[y_0,y_1]^{(b)}} 
\]
as in Theorem~(4.3.1) in \cite{GW1}.
\end{rm}
\end{Example}

\vspace{5mm}

\noindent
\begin{tabular}{l}
Mitsuyasu Hashimoto \\
Graduate School of Mathematics \\
Nagoya University \\
Chikusa-ku, Nagoya 464-8602, Japan \\
{\tt hasimoto@math.nagoya-u.ac.jp} \\
{\tt http://www.math.nagoya-u.ac.jp/\~{}hasimoto}
\end{tabular}

\vspace{2mm}

\noindent
\begin{tabular}{l}
Kazuhiko Kurano \\
Department of Mathematics \\
School of Science and Technology \\
Meiji University \\
Higashimita 1-1-1, Tama-ku \\
Kawasaki 214-8571, Japan \\
{\tt kurano@math.meiji.ac.jp} \\
{\tt http://www.math.meiji.ac.jp/\~{}kurano}
\end{tabular}

\end{document}